\documentclass[11pt,a4paper]{article}
\usepackage{amsmath}
\usepackage{bbding}

\usepackage{amsmath,amsfonts,amssymb,amsthm}

\newcommand{\diam}{\operatorname{diam}}

\newcommand{\counte}{equation}%%%设置定理与公式公用计数器
\newtheorem{theorem}[\counte]{\bf Theorem}%[section]

\newtheorem{lemma}[\counte]{\bf Lemma}

\newtheorem{coro}[\counte]{\bf Corollary}

\newtheorem{remark}[\counte]{\bf Remark}

\allowdisplaybreaks

\numberwithin{equation}{section}%%公式随section重新计数
\renewcommand{\thefootnote}{\fnsymbol{footnote}}

\begin{document}

\renewcommand{\thefootnote}{\arabic{footnote}}

\centerline{\bf\Large On the Blaschke's Conjecture
\footnote{Supported by NSFC 11471039. \hfill{$\,$}}}

\vskip5mm

\centerline{Xiaole Su, Hongwei Sun, Yusheng Wang\footnote{The
corresponding author (E-mail: wyusheng@bnu.edu.cn). \hfill{$\,$}}}

\vskip6mm

\noindent{\bf Abstract.} The Blaschke's conjecture asserts that if
$\diam(M)=\text{Inj}(M)=\frac\pi2$ (up to a rescaling) for a
complete Riemannian manifold $M$, then $M$ is isometric to $\Bbb
S^n(\frac12)$, ${\Bbb R\Bbb P}^{n}$, ${\Bbb C\Bbb P}^{n}$, ${\Bbb
H\Bbb P}^{n}$ or ${\Bbb Ca\Bbb P}^{2}$ endowed with the canonical
metric. In the paper, we prove that the conjecture is true if we in
addition assume that $\sec_M\geq1$.

\vskip1mm


\vskip1mm


\vskip6mm

\setcounter{section}{-1}

\section{Introduction}

\vskip2mm

It is well known that the sphere $\Bbb S^n(\frac12)$ and projective
spaces ${\Bbb K\Bbb P}^{n}$ endowed with canonical metrics (here the
canonical metric on a ${\Bbb K\Bbb P}^{n}$ is induced from the unit
sphere), where $\Bbb K=\Bbb R, \Bbb C, \Bbb H$ or $\Bbb
Ca$ and $n\leq 2$ if $\Bbb K=\Bbb Ca$, satisfies
$$\diam(M)=\text{Inj}(M)=\frac\pi2,\eqno{(0.1)}$$
where $\diam(M)$ and $\text{Inj}(M)$ are the diameter and injective
radius of $M$ respectively. And $\Bbb S^n(\frac12)$ and ${\Bbb K\Bbb
P}^{n}$ are the only known examples satisfying (0.1).

\vskip2mm

\noindent{\bf Blaschke's conjecture.} {\it If a complete Riemannian
manifold $M$ satisfies $(0.1)$ (up to a rescaling), then $M$ is
isometric to $\Bbb S^n(\frac12)$ or a ${\Bbb K\Bbb P}^{n}$ endowed
with the canonical metric.}

\vskip2mm

This conjecture is of long history, for which we refer to [Be], [B],
[Bo]. (Please see (1.1) below for the reason why it is called Blaschke's conjecture.)
Up to now, the conjecture is still almost open (there are only some
partial answers to it) although (0.1) is an extremely strong condition. Note that the conjecture has no
restriction to the curvature. The present paper mainly give a
positive answer to the conjecture under the additional assumption
$\sec_M\geq1$, which is stated as follows.

\vskip2mm

\noindent{\bf Main Theorem.}  {\it If a complete Riemannian manifold
$M$ satisfies $(0.1)$ and $\sec_M\geq1$, then $M$ is isometric to
$\Bbb S^n(\frac12)$ or a ${\Bbb K\Bbb P}^{n}$ endowed with the
canonical metric.}

\vskip2mm

If the curvature has upper bound, Rovenskii-Toponogov proves
that ([RT], [SSW]):

\vskip2mm

\noindent{\bf Theorem 0.1.}  {\it If a complete, simply connected
Riemannian manifold $M$ satisfies $(0.1)$ and $\sec_M\leq4$, then
$M$ is isometric to $\Bbb S^n(\frac12)$ or a ${\Bbb K\Bbb P}^{n}\
(\Bbb K\neq \Bbb R)$ endowed with the canonical metric.}

\vskip2mm

From our Main Theorem and Theorem 0.1, one can see how beautiful is
the following Berger's rigidity theorem ([CE]).

\vskip2mm

\noindent{\bf Theorem 0.2.} {\it Let $M$ be a complete, simply
connected Riemannian manifold with $1\leq\sec_M\leq 4$. If
$\diam(M)=\frac\pi2$, then $M$ is isometric to ${\Bbb
S}^{n}(\frac12)$ or a ${\Bbb K\Bbb P}^{n}\ (\Bbb K\neq \Bbb R)$
endowed with the canonical metric.}

\vskip2mm

In fact, ``$1\leq\sec_M\leq 4$'' and ``simply connected'' imply that
$\text{Inj}(M)\geq\frac\pi2$ ([CG]),
so ``$\diam(M)=\frac\pi2$'' implies that $M$ (in Theorem 0.2)
satisfies (0.1) (note that $\text{Inj}(M)\leq\diam(M)$). Hence,
the Main Theorem implies Theorem 0.2 in the premise of (0.1) (so does
Theorem 0.1). (Of course, ``$\sec_M\geq 1$'' implies that
$\diam(M)\leq\pi$. And the Maximal Diameter Theorem asserts
that {\it if $\diam(M)=\pi$, then $M$ is isometric to ${\Bbb S}^{n}(1)$},
so Theorem 0.2 is also called the Minimal
Diameter Theorem. Moreover, inspired by Theorem 0.2, Grove-Shiohama,
Gromoll-Grove and Wilhelm supply some beautiful (but not purely
isometrical) classifications under ``$\sec_M\geq 1$ and
$\diam(M)\geq\frac\pi2$ or $\text{Rad}(M)\geq\frac\pi2$'' ([GG1], [W]).)

Moreover, from the proof in [CE] for Theorem 0.2, it is not hard to see the following.

\vskip2mm

\noindent{\bf Theorem 0.2$'$.} {\it Let $M$ be a complete
Riemannian manifold satisfying $(0.1)$ and $1\leq\sec_M\leq 4$. Then $M$ is isometric to ${\Bbb
S}^{n}(\frac12)$ or a ${\Bbb K\Bbb P}^{n}$
endowed with the canonical metric.}

\vskip2mm

We will end this section with the idea of our proof of the Main
Theorem. We first prove that
$\{p\}^{=\frac\pi2}\triangleq\{q\in M||pq|=\frac\pi2\}$ for any
$p\in M$ (we denote by $|pq|$ the distance between $p$ and $q$
in the paper) is a complete totally geodesic submanifold in $M$. Then
using Theorem 1.3 below and the Toponogov's comparison
theorem, we will derive that $1\leq\sec_M\leq 4$
by the induction, and thus the proof is done by Theorem 0.2$'$.
(We would like to point out that, in the premise of
Theorem 1.3, we can use the method in [GG1-2] and [W]
to give the proof (which concerns many significant classification results).
Compared with it, our proof is much easier.)

%%%%%%%%%%%%%%%%%%%%%%%%%%%%%%%%%%%%%%% Section 1  %%%%%%%%%%%%%%%%%%%%%%%%%%%%%%%%%%%%%%%

\section{Blaschke's manifolds}

A closed Riemannian manifold $M$ is called a Blaschke's one if it
is {\it Blaschke} at each point $p\in M$, i.e. {\it $\Uparrow_q^p$
is a great sphere in $\Sigma_qM$ for any $q$ in the cut locus of
$p$} ([Be]), where $\Sigma_qM\triangleq\{v\in T_qM| |v|=1\}$
and $\Uparrow_q^p\triangleq\{\text{the unit tangent vector at
$q$ of a minimal}$ $\text{geodesic from $q$ to $p$}\}$. On a
Blaschke's manifold, one can get the following not so obvious fact
(p.137 in [Be]).

\vskip2mm

\noindent{\bf Proposition 1.1.} {\it For a Blaschke's manifold $M$,
we have that $\diam(M)=\text{\rm Inj}(M)$.}

\vskip2mm

A much more difficult observation is that (p.138 in [Be]):

\vskip2mm

\noindent{\bf Proposition 1.2.} {\it Given a closed Riemannian
manifold $M$ and a point $p\in M$, if $|pq|$ is a constant for all
$q$ in the cut locus of $p$, then $M$ is Blaschke at $p$.}

\vskip2mm

Obviously, it follows from Propositions 1.1 and 1.2 that
$$\text{\it a closed Riemannian
manifold $M$ is Blaschke\ \ $\Leftrightarrow$\ \ $\diam(M)=\text{\rm
Inj}(M)$}.\eqno{(1.1)}$$

\vskip2mm

Up to now, the Blaschke's conjecture is solved only for spheres.

\vskip2mm

\noindent{\bf Theorem 1.3 ([Be],[B]).} {\it If a Blaschke's manifold
is homeomorphic to a sphere, then it is isometric to the unit sphere (up to a rescaling).}

%%%%%%%%%%%%%%%%%%%%%%%%%%%%%%%%%%%%%%% Section 2  %%%%%%%%%%%%%%%%%%%%%%%%%%%%%%%%%%%%%%%

\section{Proof of the Main Theorem}

We first give our main tool of the paper---the Toponogov's comparison theorem.

\begin{theorem}[{[P], [GM]}]\label{2.1}
Let $M$ be a complete Riemannian manifold with $\sec_M\geq\kappa$,
and let $\Bbb S^2_\kappa$ be the complete, simply
connected $2$-manifold of curvature $\kappa$.

\noindent {\rm(i)} For any $p\in M$ and minimal geodesic $[qr]\subset M$,
we associate $\tilde p$ and a minimal geodesic $[\tilde q\tilde r]$
in $\Bbb S^2_\kappa$ with $|\tilde p\tilde q|=|pq|,|\tilde p\tilde
r|=|pr|$ and $|\tilde r\tilde q|=|rq|$. Then for any $s\in[qr]$ and
$\tilde s\in[\tilde q\tilde r]$ with $|qs|=|\tilde q\tilde s|$,
we have that $|ps|\geq|\tilde p\tilde s|$.

\vskip1mm

\noindent {\rm(ii)}  For any minimal geodesics $[qp]$ and $[qr]$ in $M$,
we associate minimal geodesics $[\tilde q\tilde p]$ and $[\tilde q\tilde r]$
in $\Bbb S^2_\kappa$ with $|\tilde q\tilde p|=|qp|$, $|\tilde q\tilde r|=|qr|$
and $\angle\tilde p\tilde q\tilde r=\angle pqr$. Then we have that
$|\tilde p\tilde r|\geq|pr|$.

\vskip1mm

\noindent {\rm(iii)} If the equality in {\rm (ii)} (or in {\rm (i)} for some $s$ in
the interior part of $[qr]$)
holds, then there exists a minimal geodesic $[pr]$ such that the triangle formed by $[qp]$, $[qr]$ and $[pr]$ bounds a
surface which is convex \footnote{We say that a subset $A$ is convex (resp. totally convex) in $M$
if, between any $x\in A$ and $y\in A$, some minimal geodesic $[xy]$ (resp. all minimal geodesics) belongs to $A$.}
and can be isometrically embedded into $\Bbb S^2_\kappa$.
\end{theorem}

In the rest of the paper, $M$ always denotes the manifold in the Main Theorem,
and $N$ denotes $\{p\}^{=\frac\pi2}\triangleq\{q\in M||pq|=\frac\pi2\}$ for an
arbitrary fixed point $p\in M$. We first give an easy
observation following from (0.1) (i.e. $\text{Inj}(M)=\diam(M)=\frac\pi2$) that
$$\text{\it For any $x\in M$,
there is a minimal geodesic $[pq]$ with $q\in N$ such that $x\in [pq]$.} \eqno{(2.1)}$$

\begin{lemma}\label{2.2}
$N$ is a complete totally geodesic submanifold in $M$;
and if $\dim(N)=0$, then $N$ consists of a single point.
\end{lemma}

\begin{remark}\label{2.3}{\rm Since $\sec_M\geq1$,
it follows from (i) of Theorem \ref{2.1} that $\{p\}^{\geq\frac\pi2}\triangleq\{q\in M||pq|\geq\frac\pi2\}$ is
totally convex in $M$. Note that $N=\{p\}^{\geq\frac\pi2}$
because $\diam(M)=\frac\pi2$, and that $N$ is closed in $M$.
On the other hand, since $M$ is a  Blaschke's manifold,
we know that $N$ is a submanifold in $M$ ([Be]).
It then follows that
$N$ is a totally geodesic submanifold in $M$.
This proof is short because we apply the proposition that
$N$ is a submanifold in $M$, which is a significant property
of a Blaschke's manifold ([Be]). Here, in order to show the importance of ``$\sec_M\geq1$'',
we will supply a proof only
based on the definition of a Blaschke's manifold.  }
\end{remark}

\noindent{\bf Proof.} From Remark \ref{2.3}, we know that
$N$ is totally convex in $M$, which implies that
$N$ consists of a single point if $\dim(N)=0$.
Hence, we can assume that $\dim(N)>0$; and for any geodesic
$\gamma(t)|_{t\in[0,\ell]}\subset N$, we need only
to show that its prolonged geodesic
$\gamma(t)|_{t\in[0,\ell+\varepsilon]}$ in
$M$ also belongs to $N$ for some small $\varepsilon>0$. Note that, without loss of
generality, we can assume that there is a unique minimal geodesic
between $\gamma(0)$ and $\gamma(\ell+\varepsilon)$.  Due to (2.1),
we can select $q\in N$ such that
$\gamma(\ell+\varepsilon)\in [pq]$. Observe that $q\neq\gamma(0)$
(otherwise, it has to hold that $\gamma(\ell)\in [pq]$ which contradicts
$\gamma(\ell)\in N$). Let $[q\gamma(0)]$ be a minimal
geodesic in $N$ (note that $N$ is
convex in $M$). By the first variation formula, it is easy to see
that
$$|\uparrow_q^{\gamma(0)}\xi|\geq\frac\pi2 \text{ (in } \Sigma_qM)  \text{ for any } \xi\in\Uparrow_q^{p}.$$
On the other hand, $\Uparrow_q^{p}$ is a great sphere in $\Sigma_qM$ because
$M$ is Blaschke at $p$ (see Proposition 1.2). It follows that
$$|\uparrow_q^{\gamma(0)}\xi|=\frac\pi2 \text{ for any } \xi\in\Uparrow_q^{p}$$
in fact. Then by (iii) of Theorem \ref{2.1}, there is a minimal geodesic
$[p\gamma(0)]$ such that the triangle formed by $[q\gamma(0)]$,
$[pq]$ and $[p\gamma(0)]$ bounds a surface
(containing $[\gamma(0)\gamma(\ell+\varepsilon)]$) which is convex and can be
isometrically embedded into $\Bbb S^2(1)$. It then has to hold
that $[\gamma(0)\gamma(\ell+\varepsilon)]=[\gamma(0)q]$ because
$[\gamma(0)\gamma(\ell)]$ belongs to $N$, and so
$[\gamma(0)\gamma(\ell+\varepsilon)]\subset N$.
\hfill$\Box$

\vskip2mm

Since $N$ is a complete totally geodesic submanifold in $M$, any minimal geodesic $[pq]$ for any $q\in N$
is perpendicular to $N$ at $q$, i.e.,
$$\Uparrow_{q}^{p}\ \subseteq (\Sigma_qN)^{=\frac\pi2} \text{ in } \Sigma_qM.\eqno{(2.2)}$$
Then from the proof of Lemma \ref{2.2}, we have the following corollary.

\begin{coro}\label{2.4}
For any minimal geodesics $[pq]$ and $[qq']\subset N$, there is a minimal geodesic
$[pq']$ such that the triangle formed by $[pq]$, $[qq']$
and $[pq']$ bounds a surface which is convex and can be
isometrically embedded into $\Bbb S^2(1)$.
\end{coro}

Moreover, the ``$\subseteq$'' in (2.2) can be changed to ``='' in fact.

\begin{lemma}\label{2.5}
For any $q\in N$, we have that $\Uparrow_{q}^{p}=(\Sigma_qN)^{=\frac\pi2}$ in $\Sigma_qM$.
\end{lemma}

\noindent{\bf Proof.} According to (2.2), it suffices to show that
for any $\zeta\in (\Sigma_qN)^{=\frac\pi2}$ there is a minimal
geodesic $[qp]$ such that $\uparrow_{q}^{p}=\zeta$. Note that there
is a minimal geodesic $[qx]$ ($x\in M$) such that
$\uparrow_q^x=\zeta$, and we can assume that there is a unique
geodesic between $q$ and $x$. And it follows from (2.1) that there
is a minimal geodesic $[pq_x]$ with $q_x\in N$ such that $x\in
[pq_x]$. Hence, we need only to show that $q_x=q.$ If this is not
true, then by Corollary \ref{2.4} there are minimal geodesics $[pq]$
and $[qq_x]\subset N$ such that the triangle formed by $[pq]$,
$[pq_x]$ and $[qq_x]$ bounds a surface $D$ which is convex and can be
isometrically embedded into $\Bbb S^2(1)$. Note
that $[qx]$ belongs to $D$. This is impossible because both $[qp]$
(see (2.2)) and $[qx]$ are perpendicular to $[qq_x]$ at $q$ (in
$D$). \hfill $\square$

\vskip2mm

Now we will give the proof of our Main Theorem.

\vskip2mm

\noindent{\bf Proof of the Main Theorem.}

Note that, according to Theorem 0.2$'$, we need only to show that
$$1\leq\sec_M\leq 4.\eqno{(2.3)}$$
We will apply the induction on $\dim(N)$.

\vskip1mm

$\dim(N)=0$: By Lemma \ref{2.2}, $N$ consists of a point, so $M$ is
homeomorphic to a sphere (because $M$ consists of minimal geodesics
between $p$ and $N$). It follows from Theorem
1.3 that $M$ is isometric to $\Bbb S^n(\frac12)$ (which implies
(2.3)).

\vskip1mm

$\dim(N)=1$: Note that $N$ is a closed geodesic of length $\pi$. Let
$q_1$ and $q_2$ be two antipodal points of $N$ (i.e. $|q_1q_2|=\frac\pi2$). It follows that
there are only two minimal geodesics between $q_1$ and $q_2$ (note
that $N$ is totally convex in $M$). Similarly, we consider
$L\triangleq\{q_2\}^{=\frac\pi2}$ containing $p$ and $q_1$, which is
a totally geodesic submanifold in $M$ of dimension $>0$ by Lemma
\ref{2.2}. Then similar to Lemma \ref{2.5}, we have that
$$\Uparrow_{p}^{q_2}=(\Sigma_pL)^{=\frac\pi2}=(\Sigma_{q_1}L)^{=\frac\pi2}
=\Uparrow_{q_1}^{q_2}.$$ This implies that there are only two
minimal geodesics between $p$ and any $q\in N$ (by Lemma \ref{2.5}).
It then is easy to see that $\sec_M\equiv1$ by Corollary \ref{2.4} (in
fact $M$ is isometric to $\Bbb{RP}^2$ with the canonical metric).

\vskip1mm

$\dim(N)>1$: Since $N$ is a complete totally geodesic submanifold in $M$
(Lemma \ref{2.2}), (0.1) implies that
$$\diam(N)=\text{Inj}(N)=\frac\pi2.\eqno{(2.4)}$$
By the inductive assumption on $N$, we have that $$1\leq\sec_N\leq
4.\eqno{(2.5)}$$ On the other hand, we claim that: {\it For any
$q\in N$, $$S(p,q)\triangleq\{\text{the point on a minimal geodesic between $p$ and $q$}\}$$
is totally geodesic in $M$ and is isometric to
$\Bbb S^m(\frac12)$, where
$m=\dim(M)-\dim(N)$}. Note that (2.3) is implied by the claim, (2.5),
Lemma \ref{2.5} and Corollary \ref{2.4}. Hence, in the rest of the
proof, we need only to verify the claim.

By (2.4), we can select $r\in N$ such that $|qr|=\frac\pi2$.
Similarly, we consider $K\triangleq\{r\}^{=\frac\pi2}$ containing
$p$ and $q$, which is a complete totally geodesic submanifold in $M$ with
$\dim(K)>0$; moreover, we have that
$$\Uparrow_{p}^{r}=(\Sigma_pK)^{=\frac\pi2},$$
and $\Uparrow_{p}^{r}$ is isometric to a unit sphere by Lemma 2.5.
On the other hand, note that $\Uparrow_{r}^{p}$ is isometric to $\Bbb S^{m-1}(1)$ by Lemma 2.5,
and that $\Uparrow_{r}^{p}$ is isometric to $\Uparrow_{p}^{r}$.
Therefore, it is easy to see (again from Lemma 2.5 on $K$) that
$$\dim(K)=\dim(N).$$
Hence, by the inductive assumption on $K$ (similar to on $N$), $K$
is isometric to $\Bbb S^l(\frac12)$ or a  $\Bbb {KP}^l$ endowed with
the canonical metric, which implies the claim above. \hfill$\Box$

%%%%%%%%%%%%%%%%%%%%%%%%%%%%%%%%%%%%%%%%%%%%%%%%%%%%%%%%%%%%%%%

\noindent School of Mathematical Sciences (and Lab. math. Com.
Sys.), Beijing Normal University, Beijing, 100875
P.R.C.\\
e-mail: suxiaole$@$bnu.edu.cn; wyusheng$@$bnu.edu.cn

\vskip2mm

\noindent Mathematics Department, Capital Normal University,
Beijing, 100037 P.R.C.\\
e-mail: hwsun$@$bnu.edu.cn

\end{document}